\documentclass[12pt]{amsart}
\usepackage{amsmath}
\usepackage{amssymb}
\usepackage{anysize}
\marginsize{1.25cm}{1.25cm}{1.75cm}{2.5cm}

\theoremstyle{plain}
\newtheorem{theor}{Theorem}

\theoremstyle{definition}

\newtheorem{claim}{Claim}

\theoremstyle{remark}

\numberwithin{equation}{section}
\newcommand{\angbr}[1]{\langle #1 \rangle}

\newcommand{\Reals}{{\mathbb R}}
\newcommand{\Rationals}{{\mathbb Q}}
\newcommand{\Poset}{{\mathbb P}}

\newcommand{\Naturals}{{\mathbb N}}

\title[Maximal Orthogonal Families]{Martin's Axiom and Maximal Orthogonal Families}
\author[S. Shelah]{Saharon Shelah}
\address{Department of Mathematics, Rutgers University, Hill Center,
 Piscataway, 
 New Jersey, U.S.A. 08854-8019}
\curraddr{Institute of Mathematics\\Hebrew University\\
Givat Ram, Jerusalem 91904, Israel}
\email{shelah@math.rutgers.edu}
\thanks{This research was supported by The Israel Science
Foundation, founded by the Israel Academy of Sciences and Humanities, and
by NSF grant No. NSF-DMS97-04477. This is number 8XX in the author's
personal numbering system.
}
\begin{document}
\begin{abstract}
It is shown that Martin's Axiom for $\sigma$-centred partial orders
implies that every maximal orthogonal family in $\Reals^\Naturals$ is of
  size $2^{\aleph_0}$. 
\end{abstract}
\maketitle

For $x,y\in\Reals^\Naturals$ define the inner product
$$\angbr{x,y}=\sum_{n=0}^\infty x(n)y(n)$$ in the obvious way
noting, however, that it may not be finite
or, indeed, may not even exist. Nevertheless, if
$\angbr{x,y}$ converges and equals $0$ then 
 $x$ and $y$ are said to be orthogonal.
A family $X\subseteq \Reals^\Naturals$ will be  said to be maximal
orthogonal
if any two of its elements are orthogonal and for every $y \in
\Reals^\Naturals\setminus X$ there is some $x \in X$ which is not
orthogonal to $y$. In \cite{step.43} various results are established
which indicate a similarity between maximal orthogonal familes and
maximal almost disjoint families of sets of integers.
There is a key distinction though: While no infinite,
 countable family of
subsets of the integers can be maximal almost disjoint, there are
countably infinite maximal orthogonal families.
 In \cite{step.43} the question
 of whether it is possible to construct a  maximal
orthogonal family of cardinality $\aleph_1$
without assuming any extra set theoretic axioms was
posed. The following theorem establishes
 that this is not possible:
 \begin{theor}
Martin's 
Axiom for $\sigma$-centred partial orders
implies that every uncountable, maximal orthogonal family in $\Reals^\Naturals$ is of  size $2^{\aleph_0}$.    
 \end{theor}
 \begin{proof}
   Let $X\subseteq \Reals^\Naturals$ be an  uncountable orthogonal family
   of cardinality less than $2^{\aleph_0}$. It will be shown that it
   can be extended to a larger orthogonal family. Before continuing, some
   notation and terminology will be established.  Whenever a topology
   on  $\Reals^\Naturals$ is mentioned this will refer to
    the usual product topology. Basic neighbourhoods of
   $\Reals^\Naturals$ will be taken to be sets of the form 
${\mathcal V}  =\prod_{i=0}^k (a_i,b_i)$ where the end points $a_i$
   and $b_i$ are all rational. The integer $k$ will be
   said to be the 
   length of $\mathcal V$ and will be denoted by
$l({\mathcal V})$  while $\max_{i\leq k}(b_i - a_i)$ will be
   referred to as the width of ${\mathcal V}$ and will be denoted by
$w({\mathcal V})$.

Let $\Poset$ be the
   set of all triples $p = ({\mathcal V},W,\eta)$ such that:
   \begin{itemize}
   \item ${\mathcal V}$ is a basic open subset of $\Reals^\Naturals$
\item $W$ is a finite subset of $X$
\item $\eta\in \Rationals$ and $\eta \geq w({\mathcal V})$
\item if $U$ is the set of all $x\in X\cap \mathcal V$ such that 
$|\sum_{i=0}^kw(i)x(i)| < \eta $
for any
 $k$  greater than the length of ${\mathcal V}$ and any $ w\in
  W$ 
then $|U| \geq \aleph_1$. 
   \end{itemize}
Define
${\mathcal V}(p) =  {\mathcal V}$, $W(p) = W$, $\eta(p) = \eta$ and
$U(p)=U$.
Define $p  \leq_\Poset p'$ if and only if
\begin{itemize}
\item ${\mathcal V}(p) \subseteq {\mathcal V}(p')$
\item $W(p) \supseteq W(p')$
\item $\eta(p) \leq\eta(p')$
\item and for each $t\in {\mathcal V}(p)$ and each integer $j$ such that
 $l({\mathcal V}(p'))
 < j \leq l({\mathcal V}(p))$ the inequality
 $|\sum_{i=0}^jt(i)w(i) |< \eta(p')$ holds for
for every $w\in W(p')$. 
\end{itemize}
Observe that $\Poset$ is $\sigma$-centred since, given any finite set of
 conditions ${\mathcal P}\subseteq \Poset$ 
such that ${\mathcal V}(p') = {\mathcal
  V}$ and $\eta(p) = \eta$ for each $p \in {\mathcal P}$, the triple
$({\mathcal V},\bigcup_{p\in{\mathcal P}}W(p),\eta)$ is a lower bound
for all of them.

It will be shown that the following sets are dense in $\Poset$:
\begin{itemize}
\item $A(x) = \{ p \in \Poset : x\in W(p)\}$
\item $B(x) = \{ p \in \Poset : x\notin {\mathcal V}(p)\}$
\item $C(m) = \{ p \in \Poset : \eta(p) < 1/m \}$
\item $D(m) = \{ p \in \Poset : l({\mathcal V}(p)) > m \}$
\end{itemize}
where $x\in X$ and $m\in \Naturals$.
Given that this assertion can be established, let $G\subseteq \Poset$
be a filter such that $$G\cap A(x)\cap B(x) \cap C(m) \cap D(m)
\neq \emptyset$$ for each $x\in X$ and
 $m\in\Naturals$.
Using that $G\cap C(m)\cap D(m)\neq \emptyset$ for each $m \in \Naturals$, 
let $x_G\in \Reals^\Naturals$ be the unique sequence such that
 $x_G \in
{\mathcal V}(p)$ for each $p \in G$. 
Observe that $x_G\neq x$ if $G\cap B(x) \neq \emptyset$. Hence $x_G
\notin X$.

To see that
$\angbr{x_G,x} = 0$ for each $x\in X$, let $x\in X$ and $\epsilon >0$ be
given and choose $k\in\Naturals$ 
such that $1/k < \epsilon$. Then select $p \in G\cap
A(x)\cap C(k)$. Now,  
given any $j$  greater than the length of ${\mathcal V}(p)$
use that $G\cap D(j)\neq \emptyset$
to choose $p'\in G\cap D(j)$ such that
$p' \leq_\Poset p$.
It is an immediate consequence of the definition of $\leq_\Poset$ and the
facts that $x_G\in {\mathcal V}(p')$, $x\in W(p)\subseteq W(p')$ and $l({\mathcal V}(p)) \leq j
\leq l({\mathcal V}(p'))$ that
$|\sum_{i=0}^jx_G(i)x(i)| < \eta(p) <
1/k < \epsilon $.  
Since $\epsilon$ was arbitrary, it follows that $\angbr{x_G,x} = 0$.

So all that remains to be shown is that the sets
$ A(x)$, $ B(x) $, $ C(m) $ and $ D(m)$ are dense
 for each $x\in X$ and
 $m\in\Naturals$. 
 \begin{claim}
$C(m)\cap D(m)$ is dense for any $m \in\Naturals$. Moreover, for any
$p\in\Poset$ and any uncountable $Z\subseteq U(p)$ it is possible to
find
$q\leq p$ in $C(m)\cap D(m)$ such that $Z\cap U(q)$ is uncountable.
 \end{claim}
 \begin{proof}
Let $p\in \Poset$ and $Z \subseteq U(p)$ be uncountable.
For each $x \in Z$ there is some $k(x)\geq m$ such that 
$|\sum_{i=0}^{j}w(i)x(i)| < 1/m $ for each $j\geq k(x)$ and $w\in W(p)$.
Choose $k$ such that $U = \{ x\in Z : k(x) = k\}$ is uncountable.
Since $\Reals^\omega$ has a countable base it is possible to find $x\in
U$ which is a complete accumulation point of $U$.
By the definition of $x\in U(p)$ it follows that
$|\sum_{i=0}^kw(i)x(i)| < \eta(p) $ for every $w \in W(p)$. 
Therefore there is some $\delta > 0$ such that
for any sequence $\{t_j\}_{j=0}^k$ such that
 $|x(j) - t_j| < \delta$ for
 each $j \leq k$ the inequality
$|\sum_{i=0}^kw(i)t_i| < \eta(p) $ holds for every $w \in W(p)$.  
 
Let ${\mathcal W}$ be a neighbourhood of $x$ with length $k$
 but of width less than the minimum of $\delta$ and $1/m$.
 Let $q = ({\mathcal W},W(p),1/m)$ 
and note that $U\cap {\mathcal W}\subseteq U(q)\cap Z$ 
and $U\cap {\mathcal W}$ is uncountable 
 since $x $ was
chosen to be a complete accumulation point of $U$.
Hence $q \in \Poset$ is as required. It is also easily
verified that  the choice of $\delta$ guarantees that
$ q\leq_\Poset p$ and that $q\in C(m)\cap D(k) \subseteq C(m) \cap D(m)$.
    \end{proof}
 \begin{claim}
$A(x)$ is dense for any $x \in X$.    
 \end{claim}
 \begin{proof}
   Let $p\in \Poset$. Choose some integer $m\geq l({\mathcal V}(p))$
   such that if $Z$ is defined to be the set of all $z\in U(p)$ such
   that 
$|\sum_{i=0}^jz(i)x(i)| < \eta(p) $ 
for each $ j\geq m$ then $|Z| \geq \aleph_1$.  Use the claim about the 
density of $C(m) \cap D(m)$ to find $ q\leq p$ such that $Z\cap U(q)$
   is uncountable and $l({\mathcal V}(q)) \geq m$.
It follows that there are uncountably many $z\in X\cap{\mathcal V}(q)$
   such that 
$|\sum_{i=0}^jz(i)x(i)| < \eta(p) $ for each $j\geq
l({\mathcal V}(q)) \geq m$.
This, in conjunction with the fact that $p\in \Poset$,  implies that 
$|\sum_{i=0}^jz(i)w(i)| < \eta(p) $ for each $j\geq
l({\mathcal V}(q))$ and $w \in W(p) \cup\{x\}$.
Therefore, if $q'$  is defined to be  $({\mathcal V}(q), W(p)
   \cup\{x\}, \eta(p))$ then
$q'\in \Poset\cap A(x)$ and $ q' \leq_\Poset p$.
 \end{proof}
 \begin{claim}
$B(x)$ is dense for any $x \in X$.    
 \end{claim}
 \begin{proof}
   Let $p\in \Poset$. For each $z\in  U(p)\setminus \{x\}$ choose a pair of integers
   $(m(z),e(z))$ such that $$|x(m(z)) - z(m(z))| > 1/e(z)$$ and then
   let $(m,e)$ be some pair of integers such that $|\{z\in U(p) :
   (m(z),e(z)) = (m,e)\}| \geq \aleph_1$. Let $k$ be the maximum of $m$
   and $e$. It follows that
    for each $z\in Z$ no
   neighbourhood ${\mathcal W}$ of $z$ of length $k$ and width $1/k$
   contains   $x$.
Use   the claim about the 
density of $C(k) \cap D(k)$ to find $ q\leq p$ such that $Z\cap
   U(q)\neq \emptyset$ and $l({\mathcal V}(q)) \geq k$. It follows $x\notin
   {\mathcal V}(q)$ and so $q \in B(x)$. 
 \end{proof}
This concludes the proofs of the claims and, hence, the proof of the theorem. \end{proof}

%\bibliographystyle{plain}
%\bibliography{../myabbrev,../steprans}

\end{document}